\documentclass[12pt]{article}
\usepackage[cp1251]{inputenc}

\usepackage{amsfonts}
\usepackage{amsmath,amssymb,amscd,latexsym,euscript,mathrsfs}
\usepackage[all]{xy}

\newcommand{\fU}{{\mathfrak{U}}}

\newcommand{\Set}{{\rm Set}}
\newcommand{\Cat}{{\rm Cat}}
\newcommand{\Ab}{{\rm Ab}}
\newcommand{\SC}{{\rm SC}}

\newcommand{\NN}{{\mathbb N}}

\newcommand{\ZZ}{{\mathbb Z}}

\newcommand{\mC}{{\mathscr C}}
\newcommand{\mD}{{\mathscr D}}
\newcommand{\mA}{{\mathscr A}}

\newcommand{\mF}{{\mathscr F}}

\newcommand{\mK}{{\mathscr K}}

\newtheorem{theorem}{\bf Theorem}[section]
\newtheorem{lemma}[theorem]{\bf Lemma}
\newtheorem{proposition}[theorem]{\bf Proposition}
\newtheorem{corollary}[theorem]{\bf Corollary}
\newtheorem{definition}[theorem]{\bf Definition}
\newtheorem{example}[theorem]{\bf Example}
\newtheorem{remark}[theorem]{\bf Remark}

\DeclareMathOperator{\che}{\check{H}}
\DeclareMathOperator{\Imm}{Im}
\DeclareMathOperator{\Ker}{Ker}

\DeclareMathOperator{\dirlim}{\underrightarrow{\rm lim}}
\DeclareMathOperator{\invlim}{\underleftarrow{\rm lim}}
\DeclareMathOperator{\pt}{pt}
\DeclareMathOperator{\Ob}{Ob}

\def\leq{\leqslant}
\def\geq{\geqslant}

\title
{
\v{C}ech cohomology of partially ordered sets
{
\normalsize{(Dedicated to the 90th anniversary of V I Kuzminov)}
}}

\author{Ahmet A. Husainov}
\date{} 

\begin{document}

\maketitle

\begin{abstract}
The article is devoted to a comparison of the \v{C}ech cohomology with the coefficients in a presheaf of Abelian groups and the topos cohomology of the sheaf generated by this presheaf for a poset with the Aleksandrov topology. The article consists of three parts.
The first part provides information from the theory of cohomology of small categories and cohomology of simplicial sets with systems of coefficients. The second part is devoted to Laudal's Theorem stating that covering cohomology for an arbitrary topological space with coefficients in the presheaf of Abelian groups is isomorphic to the derived limit functors.
The third part presents the main results.
The criterion for the invariance of cohomology groups of small categories when passing to the inverse image leads to necessary and sufficient isomorphism conditions for the \v{C}ech cohomology of an arbitrary presheaf and the topos cohomology of the sheaf generated by this presheaf.
In particular, for a finite poset, these conditions reduce to the acyclicity of the upper secrions of Dedekind-MacNeille cuts having a non-empty lower section, and the verification of these conditions is algorithmically computable.

\end{abstract}

2020 Mathematics Subject Classification 55U10, 55N05, 55N30, 18B35, 18N50, 55U05.


\tableofcontents

\section{Introduction}

The cohomology of partially ordered sets (shortly posets) with coefficients in diagrams can be considered as the cohomology of $T_0$-Aleksandrov spaces with coefficients in sheaves.
They have found various applications for generalizing homology theories for various classes of topological spaces - spectral homology theory, Borel-Moore homology theory for compact spaces, coinciding in the metric case with Steenrod-Sitnikov homology, normal Alexander-Spanier homology theory for general topological spaces \cite{deh1962}, \cite{lau1965}, \cite{kuz1967}, \cite{bal2022}.
In the article \cite{san2020} a number of theorems are proved for the homology and cohomology of co-sheaves and sheaves on finite posets with the Aleksandrov topology, previously known for locally compact topological spaces.

In our time, new problems have arisen related to applications of cohomology of posets.
For applications, it is important to have methods for calculating homological invariants.
For example, to calculate minimal sensor networks, the Euler characteristic of subsets of a finite poset \cite{tan2017} is calculated.
Sensor networks are also motivated by the problem of covering a region in Euclidean space with balls of a fixed radius in unknown locations.
The paper \cite[Theorem 4.3]{sil2007} gives a homology criterion that guarantees that a set of balls covers a bounded region. Homologies are calculated using the Vietoris-Rips complex. The problem is solved within the framework of the theory of persistent homology.
Aleksandrov's topology on a poset also found applications in persistent homology theory \cite[Chapter 5]{ber2020}; this topology is closely related to the $\gamma$-topology \cite{kas2018} on finite-dimensional Euclidean spaces.

Besides {\v{C}ech} cohomology, there are other ways to calculate homological invariants.
In particular, Kovalevsky \cite{kov1989} introduced cellular complexes as posets with a dimension function to solve various problems of image analysis and computer graphics.
He used this model \cite{kov1989} when creating image processing software in an algorithm for converting a bitmap into a list of cells.
In the article \cite[Definition 3.1]{eve2015}, Kovalevsky's definition was clarified, which formalized the property of cells for elements of a poset, and also cell cohomology with coefficients in the diagram of Abelian groups was introduced and the theorem was proved that in the case of a locally finite poset
 its cohomology sets are isomorphic to the cohomology of the cellular complex \cite[Theorem 1]{eve2015}.
  This model is the subject of the work \cite[Theorem 4.6]{kis2022}, in which the homology of the cellular complex for the regular CW complex was studied.
Using the covering theorem with homologically trivial intersections for {\v{C}ech} homology, they obtained formulas for the Borel-Moore homology of cellular co-sheaves with coefficients in the co-sheaf on the simplicial complex and on the regular CW-complex.

We are considering the following problem, posed to me by my supervisor V.I. Kuz'minov in 1977.
Let a finite topological space be given (for example, using a preorder relation matrix).
Answer the question whether {\v{C}ech} cohomology with coefficients in an arbitrary presheaf of Abelian groups will be isomorphic to the topos cohomology of the sheaf generated by this presheaf.
By topos cohomology (or Grothendieck cohomology) we mean the derived functors of a global section of a sheaf.

The following cases of isomorphism between {\v{C}ech} cohomology and topos cohomology for Aleksandrov $T_0$-spaces are known.
As pointed out by Jensen \cite[Pages 4-5]{jen1972}, in the case of a directed set, {\v{C}ech} cohomology is isomorphic to topos cohomology.
These cohomologies are also isomorphic if the poset is a lower semilattice \cite[Corollary 6.32]{ben2020}.
Note that in this case the abstract {\v{C}ech} cohomology theory, described in the review \cite{mir2011}, is developed.
If the {\v{C}ech} cohomology is considered for a sheaf, then this isomorphism will exist for any topological space in dimensions $0$ and $1$ \cite[III.3.8, Corollary 3]{gro1957} (but we
are studying the {\v{C}ech} cohomology for a presheaf, where this may not be true, as Example \ref{contra01} shows).

We give a new proof of Laudal's Theorem \cite[Page 262]{lau1965} that the cohomology of a covering with coefficients in the presheaf is isomorphic to the derived limit functors (Theorem \ref{cheaslim}).
Our method allows us to use the technique of comparing derived limit functors. For posets, this technique leads to a criterion for the isomorphism of the {\v{C}ech} cohomology with coefficients in the presheaf and the topos cohomology of the sheaf generated by this presheaf.

 For a finite poset, checking the conditions of the criterion we obtained reduces to the acyclicity of the upper sections of Dedekind-MacNeille cuts \cite[IV, \S 11, (5)]{kur1968} having a non-empty lower section, which are constructed preserving the upper and lower bounds when embedding a poset into a complete lattice.
This shows that there is an algorithm that allows, for an arbitrary given finite poset, to determine whether, for any presheaf, the canonical homomorphism between the {\v{C}ech} cohomology groups with coefficients in the presheaf and the topos cohomology groups in the sheaf generated by this presheaf is an isomorphism.
Since this algorithm contains the calculation of integer homology of finite subsets, it will not be polynomial, but its software implementation is possible.

\section{Preliminaries}

Let us write down the initial definitions and notation.

\begin{itemize}
\item $\Set$ - category of sets and mappings.
\item For any set $S$, $|S|$ denotes its cardinality.
\item For any (locally small) category $\mA$, we denote the set of morphisms $a\to b$ between $a, b\in \Ob\mA$ by $\mA(a,b)$, and the morphism functor - through $\mA(-, =): \mA^{op}\times \mA\to \Set$. If $\mA$ is an Abelian category, then by $Hom_{\mA}(a,b)$ or $Hom(a,b)$ we denote the Abelian group of morphisms.
\item
$\mA^{\mC}$ is a category of functors from the small category $\mC$ to an arbitrary category $\mA$.
Functors $F$ from a small category $\mC$ to an arbitrary $\mA$ are called object diagrams of the category $\mA$ over $\mC$ and can be denoted as the family $\{F(c)\}_{c\in \mC}$.
\item $\Delta_{\mC}A$ is a functor $\mC\to \mA$ that takes constant values $A\in \Ob\mA$ on objects, and values $1_A$ on morphisms.
\item $\NN$ is the set of non-negative integers.
\item $\Cat$ is the category of small categories and functors.
\item $0$ is a group consisting of one element.
\item $\Delta$ is the category of finite linearly ordered sets $[n]= \{0, 1, \cdots, n\}$, $n\geq 0$, and non-decreasing mappings. The category $\Delta$ is generated by morphisms of the following form:
\begin{enumerate}
\item $\partial^i_n: [n-1]\to [n]$ (for $0\leq i\leq n$) is an increasing mapping whose image does not contain $i$,
\item
$\sigma^i_n: [n+1]\to [n]$ (for $0\leq i\leq n$) is a non-decreasing surjection taking the value $i$ twice.
\end{enumerate}
\end{itemize}

Let $\Phi: \mC\to \mD$ be a functor between small categories.

For each $d\in \Ob\mD$ the left fibre \cite[Appendix 2, \S3.5]{gab1967} or comma category $\Phi$ over $d$ \cite{mac2004} is the category $\Phi/ d$, whose objects are the pairs $(c\in \Ob\mC, \beta\in \mD(S(c), d))$, and the morphisms $(c, \beta)\xrightarrow{\alpha} ( c', \beta')$ are given by $\alpha\in \mC(c,c')$ satisfying the relation $\beta'\circ S(\alpha)= \beta$.
  The left fibre has a forgetful functor $Q_d: \Phi/ d\to \mC$, acting on objects as $(c, \beta)\mapsto c$, and on morphisms as $((c, \beta)\xrightarrow{\alpha} (c', \beta')) \mapsto \alpha$.

   The right fibre or comma category $d/ \Phi$ has objects $(c\in \Ob\mC, \beta\in \mD(d, S(c)))$ and morphisms $(c, \beta) \xrightarrow{\alpha}(c',\beta')$ are defined using morphisms $\alpha\in \mC(c,c')$ such that $S(\alpha)\circ\beta = \beta '$. Forgetful functor $Q_d: d/ \Phi\to \mC$. It matches objects $(c,\beta)\mapsto c$, and morphisms
   $((c, \beta)\xrightarrow{\alpha} (c',\beta'))\mapsto \alpha$.

For an arbitrary category $\mA$, the inverse image functor $\Phi^*:  \mA^{\mD}\to \mA^{\mC}$ is defined, which maps each diagram $F\in \mA^{\mD}$ into the composition $F\Phi= F\circ \Phi\in \mA^{\mC}$, and every natural transformation $\eta: F\to F'$ into a natural transformation $\eta\Phi: F\Phi \to F'\Phi$ defined by the formula $(\eta\Phi)_{c}= \eta_{\Phi(c)}$, for all $c\in \Ob\mC$.

If $\mA$ is a cocomplete category, then the functor $\Phi^*$ has a left adjoint functor $Lan^{\Phi}: \mA^{\mC}\to \mA^{\mD}$, which is called the left Kan extension and we will define $F\in \mA^{\mC}$ on objects as $Lan^{\Phi}F (d)= \dirlim^{\Phi/ d}FQ_d$ \cite[\S10. 3 (10)]{mac2004}, and on morphisms using the universality property of the colimit functor.

If $\mA$ is a complete category, then the functor $\Phi^*$ has a right adjoint functor $Ran_{\Phi}: \mA^{\mC}\to \mA^{\mD}$ called a right Kan extension and defined by the formula $Ran_{\Phi}F(d)= \invlim_{d/\Phi}FQ_d$.

A simplicial set is a functor $X: \Delta^{op}\to \Set$.
A simplicial mapping $X\to Y$ between simplicial sets is a natural transformation.
The category of simplicial sets is denoted by $\Set^{\Delta^{op}}$.

The homology groups $H_n(X)$ of a simplicial set $X$ are defined by the formula $H_n(X)= \Ker d_n/ \Imm d_{n+1}$ as the homology of the chain complex $$ 0 \leftarrow C_0(X) \xleftarrow {d_1} C_1(X) \xleftarrow{d_2} C_2(X) \xleftarrow{d_3}\cdots $$ consisting of free Abelian groups $C_n(X)= \ZZ X_n$ with bases $X_n$ and homomorphisms $d_0= 0$ and $d_n: C_n(X)\to C_{n-1}(X)$ for all $n\geq 1$, defined on the basis elements $x\in X_n$ as $d_n(x) = \sum \limits^n_{i=0}(-1)^i d^n_i(x)$.
Here $d^n_i= X(\partial^i_n)$, 
The homology groups $H_n(X)$ of a simplicial set $X$ are defined by the formula $H_n(X)= \Ker d_n/ \Imm d_{n+1}$ as the homology of the chain complex
$$
0 \leftarrow C_0(X) \xleftarrow {d_1} C_1(X) \xleftarrow{d_2} C_2(X) \xleftarrow{d_3}\cdots $$
 consisting of free Abelian groups $C_n(X)= \ZZ X_n$ with bases $X_n$ and homomorphisms $d_0= 0$ and $d_n: C_n(X)\to C_{n-1}(X)$ for all $n\geq 1$, defined on the basis elements $x\in X_n$ as $d_n(x) = \sum \limits^n_{i=0}(-1)^i d^n_i(x)$.
Here $d^n_i= X(\partial^i_n)$.

By Eilenberg's theorem \cite[Appendix 2, Theorem 1.1]{gab1967}, the homology groups of a simplicial set are isomorphic to the homology groups of its geometric realization.

Let $\mC$ be a small category. For $n\in \NN$ {\em by length} $n$ in $\mC$ is called a finite sequence of morphisms $c_0\stackrel{\alpha_1}\to c_1\to \cdots \to c_{n-1}
\stackrel{\alpha_n}\to c_n$ in category $\mC$. For an arbitrary path of length $n\geq 1$ and an integer $i$ from the interval $0\leq i\leq n$
let us denote by $c_0 \stackrel{\alpha_1}\to c_1\to \cdots \to \widehat{c_i} \to \cdots \to c_{n-1} \stackrel{\alpha_n}\to c_n$ the path equal
$$
\left\{
\begin{array}{ll}
c_1 \stackrel{\alpha_2}\to c_2\to \cdots \to c_{n-1}
  \stackrel{\alpha_n}\to c_n~, & \mbox{ if } i=0,\\
c_0 \stackrel{\alpha_1}\to c_1\to \cdots \to
c_{i-1} \stackrel{\alpha_{i+1}\circ\alpha_i}\longrightarrow c_{i+1}
\to  \cdots \to c_{n-1}
  \stackrel{\alpha_n}\to c_n~, & \mbox{ for } 1\leq i\leq n-1,\\
c_0 \stackrel{\alpha_1}\to c_1\to \cdots \to c_{n-2}
  \stackrel{\alpha_{n-1}}\to c_{n-1}~, & \mbox{ if } i=n.
\end{array}
\right.
$$

The nerve or classifying space of a small category $\mC$ is the simplicial set $B\mC: \Delta^{op} \to \Set$ equal to the restriction of the functor $\Cat(-, \mC)$ to the subcategory $\Delta\subset \ Cat$.

The nerve of a category can be defined as a sequence of sets of $n$-simplices $B_n\mC= \Cat([n], \mC)$, boundary operators $d^i_n= \Cat(\partial^i_n, \mC): B_n\mC\to B_{n-1}\mC$, $0\leq i\leq n$, and degeneracy operators
$s^i_n= \Cat(\sigma^i_n, \mC): B_n\mC\to B_{n+1}\mC$, $0\leq i\leq n$.

Since $n$-simplices $\sigma: [n]\to \mC$ are functors, they can be defined as paths $c_0\stackrel{\alpha_1}\to c_1\to \cdots \to c_{n-1} \stackrel{\alpha_n}\to c_n$, consisting of objects $c_i= \sigma(i)$ for $0\leq i\leq n$ and morphisms $\alpha_i= \sigma(i-1<i):\sigma (i-1)\to \sigma(i)$ for $1\leq i\leq n$.

The boundary operators $d^i_n: B_n\mC\to B_{n-1}\mC$, $0\leq i\leq n$, will act as
$$
d^i_n(c_0 \stackrel{\alpha_1}\to \ldots \stackrel{\alpha_n}\to c_n)=
(c_0 \stackrel{\alpha_1}\to c_1\to \cdots \to \widehat{c_i} \to  \cdots \to c_{n-1}
  \stackrel{\alpha_n}\to c_n),
$$
and the degeneracy operators
$s^i_n: B_n\mC\to B_{n+1}\mC$, $0\leq i\leq n$, -
$$
  s^i_n(c_0\xrightarrow{\alpha_1} \ldots \xrightarrow{\alpha_n} c_n) = (c_0\xrightarrow{\alpha_1} \ldots \xrightarrow{\alpha_{i}} c_i \xrightarrow{id_{c_i}} c_i\xrightarrow{\alpha_{i+1}} \ldots \xrightarrow{\alpha_n} c_n).
$$

We define the integer homology groups $H_n(\mC)$, $n\geq 0$, of a small category $\mC$ as the homology groups of its nerve $H_n(B\mC)$.
A small category $\mC$ is called acyclic if $H_n(\mC)=0$ for all $n>0$ and $H_0(\mC)=\ZZ$.

 \section{Cohomology and homology of categories}

In this section we recall the definitions of cohomology and homology of small categories.
Examples of calculating cohomology of posets are given. Formulas for derived functors of left and right Kan extensions are considered. An Oberst criterion is formulated for an inverse image functor preserving the cohomology of small categories.
The foundations of the Gabriel-Zisman homology theory for simplicial and semi-simplicial sets are described.

\subsection{Cohomology of category with coefficients in the diagram}

Let $\mC$ be a small category.
For an arbitrary family $\{A_i\}_{i\in I}$ of Abelian groups, the Cartesian product $\prod\limits_{i\in I}A_i$ will be considered as an Abelian group of functions $\varphi: I \longrightarrow \bigcup_{ i \in I} A_i~$, taking for all $i\in I$ the values $\varphi(i)\in A_i$.


For each diagram of Abelian groups $F: \mC \rightarrow \Ab$ consider the sequence
$$
0 \to C^0({\mC},F) \xrightarrow{\delta^0} C^1({\mC},F) \xrightarrow{\delta^1}
 \ldots \to C^{n}({\mC},F) \xrightarrow{\delta^n} C^{n+1}({\mC},F) \to \ldots,
$$
consisting of Abelian groups $ C^n({\mC},F) =
\prod\limits_{c_0 \stackrel{\alpha_1}\rightarrow \cdots
\stackrel{\alpha_n}\rightarrow c_n}
F(c_n)$, $n\geq 0$, and homomorphisms
$\delta^n: C^n (\mC,F) \rightarrow C^{n+1} (\mC,F)$
operating according to the formula
\begin{multline*}
  (\delta^n \varphi) (c_0 \stackrel{\alpha_1}\rightarrow \cdots
  \stackrel{\alpha_{n+1}}\rightarrow c_{n+1}) =\\
\sum_{i=0}^n (-1)^i \varphi (c_0 \stackrel{\alpha_1}\rightarrow \cdots
\stackrel{\alpha_i}\rightarrow {\hat c}_i
\stackrel{\alpha_{i+1}}\rightarrow \cdots
\stackrel{\alpha_{n+1}}\rightarrow c_{n+1}) +\\
(-1)^{n+1} F(c_n \stackrel{\alpha_{n+1}}\rightarrow
c_{n+1} ) ( \varphi(c_0 \stackrel{\alpha_1}\rightarrow \cdots
\stackrel{\alpha_n}\rightarrow c_n ) ).
\end{multline*}

Let's set $C^n(\mC,F)=0$ for $n<0$ and introduce $\delta^{-1}=0$.
For all integers $n\geq 0$ the equalities $\delta^{n+1}\delta^n=0$ hold.
We denote the resulting cochain complex by $C^*(\mC,F)$.

\begin{definition}
Abelian groups $H^n(C^*(\mC, F))= \Ker\delta^{n}/\Imm\delta^{n-1}$, $n\geq 0$, are called {\em cohomology groups of the small category $\mC$ with coefficients in the diagram $F$} and are denoted by $H^n(\mC, F)$.
\end{definition}

\begin{proposition}\label{derlim}
The cohomology groups $H^n(\mC, F)$ are naturally isomorphic to the values of the derivatives of the limit functor $\invlim^n_{\mC}F$, for all $n\geq 0$.
\end{proposition}
This is a special case of the dual statement for homology categories \cite[Appendix 2, Proposition 3.3]{gab1967}.

\begin{example}\label{pnt}
Let $\pt=[0]$ be a category consisting of a single object $0$ and an identical morphism $1_0$, which we will represent with an arrow. Consider the functor $F: \pt\rightarrow \Ab$.
Let us calculate $\invlim^n_{\pt} F$. According to Proposition \ref{derlim} it will be isomorphic to the cohomology groups of the complex $C^*(\pt, F)$.
Let us denote $F(0)$ by $A$.

The complex $C^*(\pt,F)$ will be equal to
$$
0 \rightarrow A \stackrel{0}\rightarrow A \stackrel{1_A}\rightarrow A \stackrel{0}\rightarrow
A \stackrel{1_A}\rightarrow \cdots~,
$$
whence $\invlim^0_{\pt}F=A$ and $\invlim^n_{\pt}F=0$ for $n>0$.
\end{example}

\subsection{Cohomology of categories not containing retractions}

The example \ref{pnt} shows that the complex for calculating the cohomology of a small category consists of an infinite number of non-zero Abelian groups.
For the case where the category $\mC$ is a finite poset, this problem is solved using the following proposition.

A small category $\mC$ is said to have no retractions if for any two morphisms $\alpha$ and $\beta$ from this category, the equality of the composition $\alpha\circ\beta$ to the identity morphism implies that $\alpha$ is the identity.

In particular, if $\mC$ is a poset, then $\mC$ has no retractions.

\begin{proposition}\label{limplus}
Let $\mC$ be a small category without retractions.
Then for any diagram of Abelian groups on $\mC$ the groups $H^n(\mC, F)$ for all $n\geq 0$ are isomorphic to the $n$-th cohomology subcomplex $C^*_+(\mC, F )\subseteq C^*(\mC, F)$ consisting of products
$$
C^n_+(\mC, F)= \prod\limits_{c_0 \xrightarrow[\alpha_1]{} c_1
\xrightarrow[\alpha_2]{} \cdots \xrightarrow[\alpha_n]{} c_n}
F(c_n)
$$
over all sequences of morphisms $\alpha_i$ that are not identical for $i>0$.
\end{proposition}

 The proof is given in \cite[Proposition 2.2]{X1997}.
 \\

In what follows we will work with diagrams on the categories $\mC= I^{op}$ where $(I, \leq)$ are posets. Morphisms $a>b$ of this category will be denoted by arrows without the labels $a\to b$. We will depict the Hasse diagram using arrows directed from top to bottom.
  The complex $C^*_+(I^{op}, F)$ will be equal to $\prod\limits_{c_0 \to c_1\to\cdots\to c_n} F(c_n)$, where $c_k\in I$ for all $0\leq k\leq n$ and $c_k\not= c_{k+1}$,
  for all $0\leq k<n$.

\begin{example}\label{threep}
Let $V=\{p_0, p_1, p_2\}$ be a poset consisting of three elements ordered by the relation $p_2<p_0$ and $p_2<p_1$. Let's construct a Hasse diagram for $V^{op}$:
$$
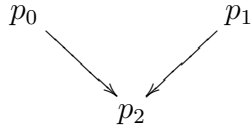
\xymatrix{
p_0\ar[rd] & & p_1\ar[ld] \\
&p_2
}
$$
For any diagram $F: V^{op}\rightarrow \Ab$ the complex $C^*_+(V^{op},F)$ consists of the Abelian groups $C^0_+(V^{op},F )=F(p_0)\times F(p_1)\times F(p_2)$, $C^1_+(V^{op},F)=F(p_2)\times F(p_2)$, and $C^n_+(V^{op},F )=0$ if $n\geq 2$.
Hence, $\invlim^n_{V^{op}} F=0$ for all $n\geq 2$.
The differential
$$
0 \rightarrow F(p_0)\times F(p_1)\times F(p_2) \stackrel{\delta^0}\rightarrow
F(p_2)\times F(p_2)\rightarrow 0
$$
operating according to formulas
\begin{itemize}
\item $(\delta^0 f)(p_0\rightarrow p_2)= f(p_2)- F(p_0\rightarrow p_2)f(p_0)$,
\item $(\delta^0 f)(p_1\rightarrow p_2)= f(p_2)- F(p_1\rightarrow p_2)f(p_0)$.
\end{itemize}

It can be represented as a matrix
$$
\left(
\begin{array}{ccc}
-F(p_0\rightarrow p_2) & 0 & 1\\
0 & -F(p_1\rightarrow p_2) & 1
\end{array}
\right)
$$
whose product by a column vector with components $f(p_i)$, $i\in \{0, 1, 2\}$, is equal to $\delta^0 f$.
Since multiplication by an invertible matrix does not change the kernel of the homomorphism, subtracting the second row from the first row, and then adding the third column multiplied by $F(p_1)\rightarrow F(p_2)$ to the second column of the resulting matrix, we arrive at the matrix
$$
\left(
\begin{array}{ccc}
-F(p_0\rightarrow p_2) & F(p_1\rightarrow p_2) & 0\\
0&0&1
\end{array}
\right)
$$
This implies that $\invlim^1_{V^{op}} F$ will be isomorphic to the quotient group
$$
F(p_2)/(\Imm F(p_0\rightarrow p_2)+ \Imm F(p_1\rightarrow p_2)~.
$$
In particular, if $F(p_0\rightarrow p_2)=F(p_1\rightarrow p_2)=0$, then $\invlim^1_V F= F(p_2)$.
\end{example}

\begin{example}\label{fourp}
Let $C_2$ be a poset defined by the following Hasse diagram for
$C^{op}_2$:
$$
\xymatrix{
	p_0 \ar[d] \ar[rd] & p_1 \ar[ld] \ar[d]\\
	p_2 & p_3
}
$$
For any functor $F: C^{op}_2\rightarrow \Ab$ the groups $\invlim^{n}_{C^{op}_2}F$ are calculated using the complex
$0\rightarrow \prod\limits_{p_i} F(p_i)\stackrel{\delta^0}\rightarrow
\prod\limits_{p_i<p_j} F(p_j) \rightarrow 0$.
The homomorphism $\delta^0$ acts from $F(p_0)\oplus F(p_1)\oplus F(p_2)\oplus F(p_3)$ to $F(p_2)\oplus F(p_2)\oplus F(p_3 )\oplus F(p_3)$ and can be defined using a matrix consisting of homomorphisms between the components of the sums
$$
\begin{array}{l}
\\
p_0\rightarrow p_2\\
p_1\rightarrow p_2\\
p_0\rightarrow p_3\\
p_1\rightarrow p_3
\end{array}
\begin{array}{cccccc}
 & p_0 & p_1 & p_2 & p_3 &\\
\hline
\vert & -F(p_0\rightarrow p_2) & 0 & 1 & 0 & \vert\\
\vert & 0 & -F(p_1\rightarrow p_2) & 1 & 0 & \vert\\
\vert & -F(p_0\rightarrow p_3) & 0 & 0 & 1 & \vert\\
\vert & 0 & -F(p_1\rightarrow p_3) & 0 & 1 & \vert
\end{array}
$$
After reduction using elementary transformations, we obtain the matrix
$$
\left(
\begin{array}{cccc}
 -F(p_0\rightarrow p_2) & F(p_1\rightarrow p_2) & 0 & 0 \\
 -F(p_0\rightarrow p_3) & F(p_1\rightarrow p_3) & 0 & 0 \\
 0 & 0 & 1 & 0 \\
 0 & 0 & 0 & 1
\end{array}
\right)
$$
Discarding the terms on which this matrix acts as an identity homomorphism, we obtain the homomorphism
$$
F(p_0)\oplus F(p_1)\rightarrow F(p_2)\oplus F(p_3)~,
$$
whose kernel is isomorphic to $\invlim_{C_2}F$, and whose cokernel is $\invlim^1_{C_2}F$.
This homomorphism acts as
$$
(a_0, a_1) \mapsto (-F(p_0\rightarrow p_2)a_0+F(p_1\rightarrow p_2)a_1, -F(p_0\rightarrow p_3)a_0+F(p_1\rightarrow p_3)a_1 )~.
$$
It follows that $\invlim^1_{C_2} F\cong \dirlim^{C_2} F$.
It is clear that $\invlim^n_{C_2} F= 0$ for all $n\geq 2$.
\end{example}

\subsection{Homology groups of small categories}

Let $\mC$ be a small category.

\begin{definition}\label{defhomolcat}
Let $F$ be a diagram of Abelian groups on $\mC$. Consider the chain complex $C_*(\mC, F)$ consisting of Abelian groups
\begin{displaymath}
 C_n({\mC},F) = \bigoplus_{c_0 \rightarrow \cdots \rightarrow c_n}
F(c_0), \quad n \geq 0,
\end{displaymath}
and homomorphisms 
$d_n= \sum\limits_{i=0}^{n}(-1)^i d^n_i:
C_n(\mC,F) \rightarrow C_{n-1}(\mC,F)$, $n>0$, 
where $d^n_i$ are defined on summand elements 
$$
(c_0 \stackrel{\alpha_1}\rightarrow c_1 \stackrel{\alpha_2}\rightarrow
\cdots \stackrel{\alpha_{n}}\rightarrow c_{n}, a)
 \in
\bigoplus_{c_0 \rightarrow \cdots \rightarrow c_n}F(c_0), \quad a\in F(c_0)
$$
by the formula
$$
d^n_i(c_0 \rightarrow \cdots \rightarrow c_n, a) =
\left\{
\begin{array}{ll}
(c_0 \rightarrow \cdots \rightarrow \widehat{c_i} \rightarrow
\cdots \rightarrow c_n, a), \quad  & $for$ ~ 1 \leq i \leq n,\\
(c_1 \rightarrow \cdots \rightarrow c_n, F(c_0\stackrel{\alpha_1}\rightarrow c_1)(a) ),
& $for$~ i = 0.
\end{array}
\right.
$$
 Homology groups or briefly homology $H_n(\mC, F)$ of the small category $\mC$ with coefficients in the diagram $F:\mC\rightarrow \Ab$ are the $n$-th homology groups of the chain complex $C_*(\mC ,F)$.
\end{definition}

\begin{proposition}
The groups $H_n(\mC, F)$ are isomorphic to the values $\dirlim^{\mC}_n F$ 
of the colimit derived functor on $F$.
\end{proposition}
The proof is given in \cite[Appendix 2, Proposition 3.3]{gab1967}.

\subsection{Derived functors of Kan extensions}

For any functor $S: \mC\rightarrow \mD$, there exist a left $Lan^S: \Ab^{\mC}\to \Ab^{\mD}$ and a right $Ran_S: \Ab^{\mC} \to \Ab^{\mD}$ Kan extensions. (The formulas for them are given in the preliminaries.)

Let us denote by $Lan^S_n$ the left $n$th derived functor of the left Kan extension and $Ran^n_S$ the right derived functor of the right Kan extensions.

\begin{lemma}\label{ranval}
Let $S: \mC\rightarrow \mD$ be a functor between small categories. For any functor $F:\mC \rightarrow \Ab$ the following formulas hold:
$$
Lan_q^S F(d)= \dirlim_q^{S/d}FQ_d~, ~ Ran^q_S F(d)=\invlim^q_{d/S}FQ_d,
\mbox{ for all } q\geq 0.
$$
\end{lemma}
The proof is given in \cite[Appendix 2, Remark 3.8]{gab1967} for $Lan^S_q$ in the case when the diagrams
over $\mC$ and $\mD$ take values in an arbitrary Abelian category $\mA$
with exact coproducs. So for $\mA=\Ab$ it is true. If we substitute $\mA= \Ab^{op}$ and take dual categories instead of $\mC$ and $\mD$ and dual to functors, we obtain the statement for $Ran^q_S$.
\hfill$\Box$

A category $\mD$ is called connected if for any of its objects $a, a'\in \Ob\mD$ there is a sequence of objects and morphisms in the category $\mD$ connecting them:
$$
	a=a_0 \to a_1 \gets a_2 \to a_3 \gets a_4\to \cdots \to a_{2k} \gets a_{2k+1}=a'.
$$

The category $\mC$ is acyclic if and only if it is non-empty and connected, and the integer homology groups $H_n(\mC)$ are equal to $0$ for all $n>0$.

\begin{theorem}\label{oberst}
Let $S: \mC\rightarrow \mD$ be a functor between small categories.
Then the canonical homomorphisms $\invlim^n_{\mD}F \rightarrow \invlim^n_{\mC}FS$ are isomorphisms for all $F\in \Ab^{\mD}$ and $n\geq 0$ if and only if for each $d\in \mD$ the comma category $S/d$ is acyclic.
\end{theorem}
In the paper \cite[Theorem 2.3]{obe1968} a dual statement was proved.
A proof using a spectral sequence is contained in \cite[Corollary 4.6]{X20232}.
A functor $S$ for which the categories $S/d$ are acyclic for all $d\in \Ob\mD$ is called strictly coinitial.

\subsection{Cohomology of (semi)simplicial sets with systems of coefficients}

Let $\Delta$ be the category of nonempty finite linearly ordered sets $[n]=\{ 0, 1, \cdots, n\}$ and nondecreasing mappings.
For an arbitrary $i\in [n]$, for $n\geq 1$, we denote by $\partial_n^i: [n-1]\to [n]$ the increasing mapping whose image does not contain $i$.
For $n\geq 0$ and $i\in [n]$ we denote by $\sigma_n^i: [n+1]\to [n]$ a non-decreasing surjection that takes twice the value $i$.
For an arbitrary diagram $F: \Delta\to \Ab$ we denote by $F[*]$ the cochain complex consisting of Abelian groups $F[n]$, where $n=0, 1, 2, \cdots$ runs through non-negative integer values. The differentials are defined as $d^n= \sum\limits_{i=0}^{n+1} (-1)^i F(\partial_{n+1}^i)$.

\begin{lemma}\label{limcomp}
For any functor $F: \Delta \rightarrow \Ab$ there are isomorphisms
$$
\invlim^n_{\Delta}F \cong H^n(F[*])
$$
\end{lemma}
This follows from the dual statement \cite[Appendix 2, Lemma 4.2]{gab1967} that the homology of the simplicial object $G: \Delta^{op}\to \mA$ of an arbitrary Abelian category $\mA$ with exact coproducts is isomorphic $\dirlim^{\Delta^{op}}_nG$. If we substitute $\mA= \Ab^{op}$, and the diagram $G= F^{op}: \Delta^{op} \to \Ab^{op}$ and go to dual categories, then $\dirlim $ will go into $\invlim$ and we will get what we are looking for.
\hfill
$\Box$
\\

A system of coefficients on a simplicial set $X$ is a functor $F: \Delta/X\rightarrow \Ab$.
Let us construct a cochain complex of Abelian groups whose cohomology is equal to $\invlim^n_{\Delta/X} F$.
To this end, consider the forgetful functor $Q: \Delta/X \to \Delta$, acting as $Q(\Delta[n]\to X)= [n]$.
Each connected component of the category $[n]/Q$ (the right layer of the functor $Q$) will have an initial object ($[n]\stackrel{1_{[n]}}\to [n], \Delta[n] \stackrel{x}\to X$).
Therefore, $Ran_Q F [n]= \prod\limits_{x\in X_n} F(x)$.

Before describing the action of $Ran_Q F$ on morphisms, let us make the following remark about the notations \cite{gab1967} and \cite{X1997}.
\begin{remark}
\cite{gab1967} uses the following notation for singular simplexes $\tilde{x}: \Delta[n]\to X$. We consider the category of singular simplices as the category of elements $x\in \coprod\limits_{n\geq 0}X_n$ of the presheaf $X$, but its morphisms are directed in the opposite direction relative to the corresponding morphisms of the category of elements.
Thus, we consider the category of singular simplices as dual to the category of elements. The morphism of singular simplices is denoted as $(\alpha,x,y)$ where $x, y\in \coprod\limits_{n\geq 0}X_n$ must satisfy $X(\alpha)(y)=x$. The corresponding morphism of the category of elements will be equal to $y\xrightarrow{\alpha}x$, where $x= X(\alpha)(y)$.
\end{remark}

To describe the action of the functor $Ran_Q F(\alpha)$ we will consider the elements of the product $\prod\limits_{x\in X_n} F(x)$ as functions of $\varphi: X_n \to \bigcup\limits_{x\in X_n } F(x)$, satisfying the condition $\varphi(x)\in F(x)$ for all $x\in X_n$.
For an arbitrary morphism $\alpha: [p]\rightarrow [q]$ of category $\Delta$ and element $z\in X_q$, the homomorphism $Ran_Q F (\alpha)$ must make the diagram commutative
\begin{equation}\label{maincosimp}
\xymatrix{
\prod\limits_{x\in X_p} F(x) \ar[d]_(.60){pr_{X(\alpha)(z)}} \ar[rr]^{(Ran_Q F) (\alpha)}
&&  \prod\limits_{y\in X_q} F(y) \ar[d]^(.60){pr_z} \\
F(X(\alpha)(z)) \ar[rr]_{\quad F(\alpha, X(\alpha)(z), z)}&& F(z)
}
\end{equation}
In other words, this homomorphism will associate with each function $\varphi\in \prod\limits_{x\in X_p} F(x)$ a function taking on $z\in X_q$ the following values
\begin{equation}\label{opsmaincosimp}
(Ran_Q F) (\alpha) (\varphi) (z) =
F(\alpha, X(\alpha)(z), z) (\varphi({X(\alpha)(z)}))
\end{equation}
In particular, for $p=n$, $q=n+1$ and $\alpha=\partial^i_{n+1}: [n]\rightarrow [n+1]$, we will have
\begin{equation}\label{ranmorph}
   (Ran_Q F) (\partial_{n+1}^i) (\varphi) (z) =
F(\partial_{n+1}^i, X(\partial_{n+1}^i)(z), z) (\varphi({X(\partial_{n+1}^i)(z) }))
\end{equation}
The functor $Ran_Q F: \Delta \to \Ab$ is equal to a cosimplicial Abelian group.

\begin{proposition}\label{ranq}
For an arbitrary simplicial set $X\in \Set^{\Delta^{op}}$ and the functor $F: \Delta/X \rightarrow \Ab$ there are isomorphisms
$$
\invlim^n_{\Delta/X} F \cong H^n(Ran_Q F[*])~,
$$
where $Ran_Q F[*]$ is a complex consisting of Abelian groups $\prod\limits_{n\in X_n} F(x)$ with differentials $d^n= \sum\limits_{i=0}^{n +1} (-1)^i Ran_Q F (\partial_{n+1}^i)$, $n\geq 0$, defined in accordance with the formulas (\ref{ranmorph}).
\end{proposition}
The proof of the dual statement in the Abelian category with exact coproducts is given in \cite[Appendix 2, Proposition 4.2]{gab1967}. \hfill $\Box$
\\

To study the cohomology of simplicial complexes, we consider the subcategory $\Delta_+\subset \Delta$, the set of objects of which coincides with the set of objects of the category $\Delta$, and the morphisms $[m]\to [n]$ are increasing mappings. Objects of the category $\Set^{\Delta^{op}_+}$ are called semisimplicial sets.

Consider the diagram of Abelian groups $F: \Delta_+\to \Ab$.
Let $F[*]$ be a cochain complex of Abelian groups $F[n]$, $n\geq 0$, and differentials $d^n= \sum\limits_{i=0}^{n+1} (- 1)^i F(\partial_{n+1}^i)$.

\begin{lemma}\label{limcompplus}
For any functor $F: \Delta_+ \rightarrow \Ab$ there are isomorphisms
$$
\invlim^n_{\Delta_+}F \cong H^n(F[*])
$$
\end{lemma}
The proof is given in \cite[Lemma 1.3]{X1997}.
\\

For a diagram $F: \Delta_+/X\to \Ab$ its right Kan extension $Ran_Q F: \Delta_+\to \Ab$ along $Q: \Delta_+/X\to \Delta_+$ will be defined in the same way as for systems of coefficients on simplicial sets.

\begin{proposition}\label{ranqplus}
For an arbitrary semisimplicial set $X\in \Set^{\Delta^{op}_+}$ and the functor $F: \Delta_+/X \rightarrow \Ab$ there are isomorphisms
$$
\invlim^n_{\Delta_+/X} F \cong H^n(Ran_Q F[*]).
$$
\end{proposition}
The proof is given in \cite[Proposition 1.4]{X1997}.

\section{Cohomology of simplicial schemes and coverings}

This part is devoted to the cohomology of a simplicial scheme and the cohomology of a covering of a topological space with coefficients in a presheaf of Abelian groups.

A complex for the cohomology of a simplicial set of singular simplices is constructed (Proposition \ref{cohconstr}) and it is proved that these cohomologies are isomorphic to the cohomology of a simplicial scheme with coefficients in the presheaf of Abelian groups (Corollary \ref{iso2}).
A new proof of Laudal's theorem is given that the cohomology of a covering of a space by open sets with coefficients in the presheaf is isomorphic to the derived functors of a limit over a category consisting of finite intersections of these open sets (Theorem \ref{cheaslim}).This allows us to construct an ordered complex for calculating the cohomology of a covering (Corollary \ref{compvalue}) reducing the amount of computation.

\subsection{Cohomology of a simplicial scheme}

A {\em simplicial scheme} or {\em simplicial complex} is a pair $(E, \mK)$ consisting of an arbitrary set of {\em vertices} $E$ and a set $\mK$ of non-empty finite subsets in $E$, satisfying the following conditions:
\begin{enumerate}
\item $(\forall{x\in E})~ \{x\}\in \mK$;
\item $(\forall S\in \mK)~ S'\not= \emptyset ~\&~ S'\subset S \Rightarrow S'\in \mK$.
\end{enumerate}
Sets $S\in \mK$ are called simplexes.
A non-negative integer $\dim S= |S|-1$ is called the dimension of the simplex $S$.
For $n\geq 0$,  $\mK_n$ denotes the set of simplices of dimension $n$.
  The simplicial scheme $(E, \mK)$ will admit the shorthand notation $\mK$. The set of vertices $E$ will be identified with the set of $0$-dimensional simplexes $\mK_0$.
\\

Let $\mK$ be a simplicial scheme, and let $(\mK, \subseteq)$ be the set of its simplices $S\in \mK$ partially ordered by the inclusion relation.
We will consider this poset as a small category in which the morphisms are the inclusions $S\subseteq S'$.

Let us denote by $\SC$ the category of simplicial schemes in which the morphisms $f: \mK\to \mK'$ are given by mappings of sets of vertices $f_0: \mK_0\to \mK'_0$ such that the image 
$f_0(S) := \{f_0(x) | x \in S\}$ of each $S\in \mK$ is a simplex of $\mK'$.

The standard simplicial schema $\Delta_n$ has vertices $\{0, 1, \cdots, n\}$. Its simplexes are all non-empty subsets of the set of vertices.

A morphism of simplicial schemes $\Delta_n\to \mK$ is called a singular simplex in $\mK$. Morphisms between singular simplices $\Delta_p\to \mK$ and $\Delta_q\to \mK$ are defined using nondecreasing maps of the sets $\{0, 1, \ldots, p\}\to \{0, 1, \ldots , q\}$, making the following diagram of simplicial schemes commutative
$$
\xymatrix{
\Delta_p \ar[rd] \ar[rr] && \mK\\
& \Delta_q \ar[ru]
}
$$

Consider the functor $T: \Delta \to \SC$, from the category of finite linearly ordered sets and nondecreasing mappings into the category of simplicial schemes, taking the values $T[n]= \Delta_n$ and $T(\alpha)=\alpha$.
The category $T/\mK$ will be equal to the category of singular simplices in $\mK$.

The system of coefficients on the category of singular simplices in $\mK$ is the functor $F: T/\mK\to \Ab$. Cohomology groups $H^n_T(\mK, F)$ are the right derived functors of the limit $\invlim^n_{T/\mK}F$, for all $n\geq 0$.

The functor $Q_{\mK}: T/\mK \to \Delta$, $(([n], T[n])\to \mK)\mapsto [n]$, is a discrete Grothendieck 
fibration, in the sense that the embedding $Q^{-1}_{\mK}([n])\subseteq [n]/Q_{\mK}$ has a right adjoint.
This allows us to construct a complex whose cohomology is equal to $H^n_T(\mK,F)$.
We will consider a construction that is even simpler.
The next lemma shows that the cohomology $H^n_T(\mK, F)$ is isomorphic to the Gabriel-Zisman cohomology of the simplicial set corresponding to the simplicial scheme $\mK$.

\begin{lemma}\label{catsingsimp}
The category $T/\mK$ is isomorphic to the category $\Delta/R\mK$, where $R\mK: \Delta^{op}\to \Set$ is a simplicial set consisting of the sets $R\mK_n= \SC( T[n], \mK)$ and mappings $R\mK(\alpha): \SC(T[q], \mK)\to \SC(T[p], \mK)$ defined for $\alpha : [p]\to [q]$ as $(y: [q]\to \mK_0) \mapsto y \circ T(\alpha)$.
\end{lemma}

Consider the functor $\Imm: T/\mK\to (\mK, \subseteq)$ that associates each singular simplex $x: T[n] \to \mK$ with its image $\Imm x= x[n]\in \mK$.

According to Lemma \ref{catsingsimp}, to construct a cohomology complex for the category $T/\mK$, it is enough to construct it for $\Delta/R\mK$.
We will need a complex for calculating the cohomology groups $\invlim^n_{\Delta/R\mK} (F\circ\Imm)$ of the category of singular simplices of a simplicial scheme.
We build this complex using the method shown in the diagram (\ref{maincosimp}) and using the formulas (\ref{opsmaincosimp}) and (\ref{ranmorph}).

We consider the simplicial set $X= R\mK$. The functor $Q$ is equal to $Q_{R\mK}: 
\Delta/R\mK\to \Delta$.
Cohomology with coefficients in $F\Imm: \Delta/R\mK\to \Ab$ is considered.

The complex consists of Abelian groups $\prod\limits_{x\in R\mK_n}F\Imm(x)$, 
where the simplices $x: \Delta_n\to \mK$ can be defined as tuples of vertices of the simplicial 
scheme $\mK$, equal to $(x_0, \ldots, x_n)$, consisting of vertices $x_i= x(i)$. The functor $\Imm: \Delta/R\mK\to \mK$ associates with every singular simplex $(x_0, \ldots, x_n)$ the simplex $\{x_0, \ldots, x_n\}\in \mK$ admitting identified vertices ($\Imm$ does not preserve dimension!).
Coboundary operators are constructed as making the following diagrams commutative,
for all $0\leq i\leq n+1$ and $z\in R\mK_{n+1}$,
where $x:=(x_0, \ldots, x_{n})$,
$y:=(y_0, \ldots, y_{n+1})$, $z:=(z_0, \ldots, z_{n+1})$,
$d^n_i:={Ran_{Q_{R\mK}}
  (F\Imm) (\partial^i_{n+1})}$.
\begin{equation}\label{maincosimp2}
\xymatrix{
\prod\limits_{x \in R\mK_n} F\Imm(x)
 \ar[d]^(.60){pr_{R\mK(\partial^i_{n+1})(z)}}
 \ar[rrrr]^{d^n_i}
&&&&  \prod\limits_{y\in R\mK_{n+1}}
F\Imm(y)
\ar[d]^(.60){pr_{z}} \\
F\Imm(R\mK(\partial^i_{n+1})(z))
 \ar[rrrr]_{\quad F\Imm(\partial^i_{n+1},
 R\mK(\partial^i_{n+1})(z), z)}
&&&& F\Imm(z)
}
\end{equation}

The bottom line will be equal to
$$
F\{z_0, \ldots, \hat{z_i}, \ldots, z_{n+1}\}
\xrightarrow{F(\{z_0, \ldots, \hat{z_i}, \ldots, z_{n+1}\}\subseteq \{z_0, \ldots, z_{n+1}\})}
F\{z_0, \ldots, z_{n+1}\}
$$

We obtain the following statement about the complex for calculating the cohomology groups of the category of singular simplexes

\begin{proposition}\label{cohconstr}
The cohomology groups $\invlim^n_{\Delta/R\mK}F\Imm$ of the category of singular simplices of the simplicial scheme $\mK$ with coefficients in $F$ for all $n\geq 0$ are isomorphic to the cohomology groups of the complex consisting of Abelian groups
$$
C^n(R\mK, F\Imm)= \prod\limits_{\{(x_0, \ldots, x_n)|\{x_0, \ldots, x_n\}\in \mK_n\}}
F\{x_0, \ldots, x_n\}
$$
and differentials
$d^n=\sum\limits_{0\leq i\leq n+1}(-1)^i d^n_i$
obtained from coboundary operators 
$d^n_i: C^n_T(\mK, F\Imm)\to C^{n+1}_T(\mK, F\Imm)$ 
defined as
\begin{multline*}
 d^n_i (\varphi) (z_0, \ldots, z_{n+1}) =
 F(\{z_0, \ldots, \hat{z_i}, \ldots, z_{n+1}\}\subseteq \{z_0, \ldots, z_{n+1}\})\\ (\varphi(z_0, \ldots, \hat{z_i}, \ldots, z_{n+1}))
\end{multline*}
\end{proposition}
\hfill$\Box$

Let us prove that the functor $\Imm: T/\mK\to (\mK, \subseteq)$ is strictly coinitial.

\begin{proposition}\label{coinsingtosc}
The left fibers of $\Imm/S$ are contractible for all simplices $S\in \mK$.
\end{proposition}
{\sc Proof.} For an arbitrary set $E$, consider the anti-discrete category $E^a$. Its set of objects is equal to $E$, and for any $x,y\in E$ the set of morphisms $E^a(x,y)$ consists of a single element. Its nerve is contractible, and therefore the nerve of category $\Delta/ E^a$ is contractible.

The set $S$ is a simplex of $\mK$, which means it is finite.
The category $\Imm/S$ will be isomorphic to the category $T/S_*$, where $S_*$ is a simplicial scheme with a set of vertices $S$ whose simplices are all non-empty subsets of $S$.

The categories $T/S_*$ and $\Delta/S^a$ are isomorphic, and therefore the nerve of the category $T/S_*$ is contractible.
Consequently, the nerve of the category $\Imm/S$ is contractible.
\hfill$\Box$
\\

A system of coefficients on a simplicial circuit $\mK$ is an arbitrary functor $F: (\mK, \subseteq)\to \Ab$ \cite[I.3.3]{god1958}.

Since there is a canonical isomorphism between the categories $T/\mK$ and $\Delta/R\mK$, we will identify these categories.

Using Theorem \ref{oberst} and Proposition \ref{coinsingtosc} we arrive at the following isomorphism, which could also be obtained using homotopy equivalence, similar to \cite[Theorem 3.1]{ben2020}.
\begin{corollary}\label{iso2}
For any simplicial scheme $\mK$ and coefficient system $F: (\mK, \subseteq)\to \Ab$ there are natural isomorphisms $\invlim^n_{(\mK, \subseteq)}F \xrightarrow{\cong} \invlim^n_{\Delta/R\mK} (F\circ\Imm)$, $n\geq 0$.
\end{corollary}

\subsection{Cohomology of open coverings}

In this subsection we present our proof of the auxiliary Theorem \ref{cheaslim}, which allows us to calculate the {\v{C}ech} cohomology for an arbitrary open covering of a topological space using derived limit functors.
This theorem was obtained by Laudal \cite[Page 262]{lau1965}.

  For an arbitrary set $X$, its cover $\fU$ is the set of subsets $U\subseteq X$ whose union $\cup\fU= \bigcup\limits_{U\in \fU}U$ is equal to $X$.

Let $X$ be an arbitrary topological space, $\fU$ its covering by non-empty open sets $U\subseteq X$. By definition, a covering is a set, which means that all elements from $\fU$ are distinct.
{\em The nerve of a covering} $\fU$ is a simplicial scheme $\mK(\fU)$ whose simplices are non-empty finite subsets $S\subseteq \fU$ such that $\cap S= \bigcap\limits_{U\in S}U \not= \emptyset$.

Let $R\mK(\fU)$ denote the simplicial set of singular simplices of the nerve of the covering $\fU$.
The category $\Delta/ R\mK(\fU)$ will be equal to the category of singular simplexes, and $(\mK(\fU), \subseteq)$ will be equal to the category of nerve simplexes.

Let $\tau_X$ be the set of all open subsets of $X$, ordered by inclusion.
An arbitrary functor $\mF: \tau_X^{op}\rightarrow \Ab$ is called a presheaf.

{\em Cohomology groups $H^n(\fU, \mF)$ of a covering with coefficients in the presheaf $\mF$} are the cohomology groups of the complex
$$
	C^n(\fU, \mF)= \prod\limits_{(U_0, \cdots, U_n)} \mF(U_0\cap \cdots \cap U_n)~,
$$
where the products are taken over all $(U_0, \cdots, U_n)$ such that $U_0\cap \cdots \cap U_n\not= \emptyset$, and the differentials $d^n: C^n(\fU, \mF)\to C^{n+1}(\fU, \mF)$
defined as
\begin{multline*}
	(d^n f)(U_0, \cdots, U_{n+1})= \\
	\sum_{k=0}^{n+1} (-1)^k \mF(U_0\cap \cdots \cap \widehat{U_k} \cap \cdots \cap U_{n+1}
		\supseteq U_0 \cap \cdots \cap U_{n+1})\\
			f(U_0, \cdots, \widehat{U_k}, \cdots, U_{n+1})
\end{multline*}
Following \cite[ch. II, $\S$5.1]{god1958}, we define the coefficient system as the functor $\mF_{\fU}:(\mK(\fU), \subseteq) \rightarrow \Ab$, setting $\mF_{\fU}( S)= \mF(\bigcap\limits_{U\in S}U)$.
From the definition of the complex $C^n(\fU, \mF)$ it follows that
$$
	C^n(\fU, \mF)= \prod\limits_{\dim\sigma=n} \mF_{\fU}(\Imm\sigma)~,
$$
where the products are taken over ${\sigma}\in \Delta/R\mK(\fU)$.

Let $\widetilde{\fU}\subseteq \tau_X$ be the set consisting of non-empty finite intersections of sets belonging to the covering $\fU$. (The set $\widetilde{\fU}\subseteq \tau_X$ consists of various subsets.)
The following functors are defined
$$
\Delta/R\mK(\fU) \stackrel{\Imm}\rightarrow (\mK(\fU), \subseteq)
 \stackrel{\overline{\cap}}\rightarrow
\widetilde{\fU}^{op}	\stackrel{\mF\vert_{\widetilde{\fU}}}\longrightarrow \Ab~.
$$
where $\overline{\cap}: (\mK({\fU}), \subseteq) \to \widetilde{\fU}^{op}$ is a functor that associates each simplex $S\in \mK(\fU )$ intersection $\cap S= \cap_{U\in S}U$.

\begin{theorem}\label{cheaslim}
Let $(X,\tau_X)$ be a topological space and $\fU$ be its covering by non-empty open subsets.
Then for any presheaf of Abelian groups $\mF$ over $X$ and $n\geq 0$ the 
following isomorphisms hold
$$
H^n(\fU,\mF) \cong \invlim^n_{\widetilde{\fU}^{op}} \mF\vert_{\widetilde{\fU}}.
$$
\end{theorem}
{\sc Proof.} First, from Proposition \ref{ranq} and Proposition \ref{cohconstr} it follows that the groups $H^n(\fU,\mF)$ will be isomorphic to $\invlim^n_{\Delta /R\mK(\fU)}(\mF_{\fU}\circ\Imm)$.
By Corollary \ref{iso2} they are isomorphic to $\invlim^n_{(\mK({\fU}), \subseteq)}\mF_{\fU}$.
The equality $\mF_{\fU}= \mF|_{\widetilde{\fU}}\circ \overline{\cap}$ is true.
This means that it remains to prove that the functor $\overline{\cap}: (\mK(\fU), \subseteq) \rightarrow \widetilde{\fU}^{op}$ is strictly coinitial.

Let $W\in \widetilde{\fU}$.
The poset $\overline{\cap}/ W$ consists of $S\in \mK(\fU)$ such that $\cap S \supseteq W$. It is ordered by the inclusion relation $S\subseteq S'$.
There are $U_0$, $\ldots$, $U_n$ for which $U_0\cap \cdots \cap U_n=W$. Therefore the category $\overline{\cap}/ W$ is non-empty.
If $\cap S\supseteq W$ and $\cap S'\supseteq W$, then $\cap (S\cup S') \supseteq W$.
This means that $\overline{\cap} / W$ is a directed set.
Therefore, $H_n(\overline{\cap} /W)\cong H_n(\pt)$ for all $n\geq 0$.
\hfill
$\Box$

\subsection{Ordered complex for cohomology of a simplicial scheme}

To calculate the cohomology groups of a simplicial scheme, we may need to introduce a linear order relation on $E$.
We need to prove that cohomology does not depend on this relation.

Let us choose this relation and denote it by the symbol $<$. Consider the semisimplicial set $\mK^{ord}: \Delta_+^{op}\to \Set$, whose $n$-dimensional simplices are the sets of elements $x_0< x_1< \cdots < x_n$ of $E$ such that $\{x_0, x_1, \ldots, x_n\}\in \mK$.
The boundary operators of this semisimplicial set are denoted by $\delta^n_i: \mK^{ord}_n\to \mK^{ord}_{n-1}$. They are defined for all $n>0$ and $0\leq i\leq n$, and remove the element $x_i$ from each such set:
$$
    \delta^n_i(x_0< x_1< \cdots <x_n) = (x_0 < \cdots <x_{i-1}< x_{i+1}< \cdots <x_n).
$$

\begin{proposition}
The category of singular simplices $\Delta_+/\mK^{ord}$ is isomorphic to the category $(\mK, \subseteq)$.
\end{proposition}
{\sc Proof.} For each singular simplex $x: \Delta_+[n]\to \mK^{ord}$ there is a bijective correspondence to the simplex $(x_0< x_1< \ldots< x_n)$ of a semisimplicial set whose elements form a simplex of $\mK$, which we denote by $\Imm x$.
The morphism between singular simplices $x\to x'$ will correspond to the inclusion $\Imm x \subseteq \Imm x'$. This correspondence is one-to-one.
We obtain the functor $\Imm: \Delta_+/\mK^{ord} \to (\mK, \subseteq)$, bijective on objects and on morphisms. This functor will be an isomorphism.
\hfill
$\Box$

Since $\Imm$ is an isomorphism, the following is true
\begin{proposition}\label{barycenter}
For an arbitrary functor $F: (\mK, \subseteq)\rightarrow \Ab$ there are isomorphisms
\begin{equation}\label{isobarycenter}
\invlim^n_{(\mK, \subseteq)} F\cong \invlim^n_{\Delta_+/\mK^{ord}} (F\circ \Imm)~.
\end{equation}
\end{proposition}

Proposition \ref{barycenter} shows that introducing a linear order relation at the vertices of a simplicial scheme can greatly reduce the amount of computation of its cohomology. A complex of cochains whose cohomology is equal to the right side of the isomorphism (\ref{isobarycenter}) consists of Abelian groups
$$
C^n_<= \prod\limits_{x_0<\cdots< x_n}F\{x_0, \ldots, x_n\}, ~n\geq 0.
$$
Its differentials  $d^n:C^n_<\to C^{n+1}_<$ are equal
$d^n= \sum\limits^{n+1}_{i=0}(-1)^i d^n_i$ where
\begin{multline*}
d^n_i(\varphi)(z_0< \ldots< z_{n+1})= \\
F(\{z_0, \ldots, z_{i-1}, z_{i+1}, \ldots, z_{n+1}\}\subseteq \{z_0, \ldots, z_{n+1}\})
\\
 (\varphi(z_0< \ldots< z_{i-1}< z_{i+1}< \ldots< z_{n+1}))
\end{multline*}

For an arbitrary presheaf of Abelian groups $\mF$ on a topological space $X$ and a covering $\fU$ of the space $X$ consisting of various non-empty open subsets, consider the simplicial scheme $\mK(\fU)$ and the system of coefficients on it $ F(U_0, \ldots, U_n):= \mF(U_0\cap \ldots \cap U_n)$, for all simplexes $\{U_0, \ldots, U_n\}\in \mK(\fU)$.
We obtain the following statement.

\begin{corollary}\label{compvalue}
For an arbitrary linear order relation $<$ on $\fU$, the cohomology groups $H^n(\fU, \mF)$ are isomorphic to the cohomology groups of the complex consisting of Abelian groups $C^n_<(\fU, \mF)= \prod \limits_{U_0<\ldots<U_n}\mF(U_0\cap\ldots \cap U_n))$ and differentials
\begin{multline*}
d^n(\varphi)(U_0<\ldots< U_n)=
 \sum\limits^{n+1}_{i=0}(-1)^i
\mF({U_0\cap \ldots \cap \hat{U_i}\cap \ldots \cap U_n}\supseteq{U_0\cap \ldots  \cap U_n})\\
\varphi(U_0< \ldots <U_{i-1}< U_{i+1}< \ldots <U_n)
\end{multline*}
\end{corollary}

\section{{\v{C}ech} cohomology for Aleksandrov spaces}

In this section, we first recall the assertions proven in \cite{deh1962}, \cite{lau1965}, \cite{jen1972} that diagrams on posets can be considered as sheaves over the $T_0$-Aleksandrov space in the sense that the topos cohomology of sheaves on this space is isomorphic to the 
values of limit derived functors on the corresponding diagram.
After this we move on to the {\v{C}ech} cohomology of a poset. We give examples showing that in the general case {\v{C}ech} cohomology and topos cohomology are not isomorphic.
We prove the main results of the work - Theorem \ref{miso} and Corollary \ref{misof} and provide examples of application.

\subsection{Cohomology of sheaves on Aleksandrov spaces}
Let $(X, \tau_X)$ be a topological space, $\mF$ a presheaf of Abelian groups over $X$,
and $(U_{\alpha})_{\alpha\in J}$
a family of open subsets $\in X$. Denote $U= \bigcup_{\alpha\in J}U_{\alpha}$.
Consider the mappings $\mF(U)\stackrel{r_{\alpha}}\rightarrow \mF(U_{\alpha})$ associated with the embeddings $U_{\alpha}\subseteq U$ and the mappings
$$
\mF(U_{\alpha}) \stackrel{p_{\alpha\beta}}\rightarrow \mF(U_{\alpha}\cap U_{\beta}), \quad
\mF(U_{\beta}) \stackrel{q_{\alpha\beta}}\rightarrow \mF(U_{\alpha}\cap U_{\beta}),
$$
defined for all $(\alpha,\beta)\in J\times J$.
We get the diagram
\begin{equation}\label{sheaf}
{\mF(U)} \stackrel{r}\rightarrow
\prod\limits_{\alpha\in J} \mF(U_{\alpha})
\xymatrix{\ar@<1ex>[r]^p \ar@<-1ex>[r]_q &}
\prod\limits_{(\alpha,\beta)\in J\times J}
\mF(U_{\alpha}\cap U_{\beta})
\end{equation}
such that
$r(g)_{\alpha}=r_{\alpha}(g)$, for $g\in \mF(U)$, 
$p(f)_{(\alpha,\beta)}= p_{\alpha\beta}(f_\alpha)$ and
$q(f)_{(\alpha,\beta)}= q_{\alpha\beta}(f_\beta)$
for each family
$f= (f_{\gamma})_{\gamma\in J}$.
A presheaf $\mF$ is called a {\em sheaf} if $(\mF(U), r)$ is an equalizer of the pair $(p,q)$ for every family of open subsets $(U_{\alpha})_{ \alpha \in J}$.
\\

Let $(X, \tau_X)$ be a topological space.
For any sheaf $\mF: \tau^{op}_X\to \Ab$ the cohomology groups $H^n(X, \mF)$ are the right derived of the functor $\Gamma_X: Shv(X, \Ab)\to \Ab$ associating a sheaf
 $\mF$ with a group $\mF(X)$. A sheaf $\mF$ is said to be flasque or flabby if the homomorphisms $\mF(X)\to \mF(U)$ are surjective for every open $U\subseteq X$.
There are enough injective objects in $Shv(X, \Ab)$, and according to \cite[Proposition 3.3.2]{gro1957} flasque sheaves satisfy the conditions \cite[Lemma 3.3.1 and its corollary]{gro1957}, whence cohomology with coefficients in a sheaf can be calculated using the flasque resolution of the sheaf.

Let $(I,\leq)$ be a poset.
Let us denote $\Lambda_i=\{j\in I: j\leq i\}$, for all $i\in I$.
We declare a subset $U\subseteq I$ {\em open} if it, together with any of its elements $i\in U$, contains all elements $j\in I$ for which $j\leq i$.
Any open subset will be equal to the union of subsets of the form 
$\Lambda_i$.
It is easy to see that the set $\tau_I$, whose elements are all open subsets of $U\subseteq I$, turn the set $I$ into the topological space $(I,\tau_I)$.

Let $Shv(I,\Ab)$ be the category of sheaves over the topological space 
$(I,\tau_I)$.
Consider the functor $\Lambda: I\to \tau_I$, defined on objects $i\in I$ as 
$\Lambda(i)= \Lambda_i$. It maps morphisms $i\leq j$ into embedding morphisms $\Lambda_i\subseteq \Lambda_j$. We can denote $\Lambda_i$ by $\Lambda(i)$.
There is a pair of adjoint functors 
$$
{(-)\circ\Lambda^{op}}: \Ab^{\tau^{op}_I}\rightleftarrows \Ab^{I^{op}} :Ran_{\Lambda^{op}}, 
$$ 
and the functor $(-)\circ\Lambda^{op}$ is exact, and $Ran_{\Lambda^{op}}$ is a complete embedding.

Let us denote by $\hat\mF$ the sheaf generated by the presheaf $\mF$.

\begin{lemma}\label{shpresh}
The category $Shv(I,\Ab)$ is isomorphic to the category $\Ab^{I^{op}}$. The sheaf $\hat\mF$ is equal to $Ran_{\Lambda^{op}}(\mF\circ\Lambda^{op})$.
\end{lemma}
{\sc Proof.} Each functor $F: I^{op}\rightarrow \Ab$ corresponds to a presheaf $Ran_{\Lambda^{op}}F$. Let us denote it by $\tilde{F}$ and show that it will be a sheaf.
For each $U\in \tau_I$ the group $\tilde{F}(U)= \invlim_{U^{op}} F\vert_U$ consists of strands of elements of the groups $F(i)$.
Embeddings $U\supseteq V$ go into homomorphisms 
$\invlim F\vert_U \rightarrow \invlim F\vert_V$ that map each thread $(f_u)_{u\in U}$ into a thread $(f_u)_{u \in V}$.
In the diagram (\ref{sheaf}), a homomorphism $r$ will take $f$ to the family $f\vert_{U_{\alpha}}$. As in the example considered above, it will implement an isomorphism between the Abelian group $\widetilde{F}(U)$ and the equalizer of the pair $(p,q)$. Therefore $\widetilde{F}$ is a sheaf.
The inverse mapping assigns to each pencil $\mF$ a diagram $\{\mF(\Lambda_i)\}_{i\in I}$.
Consequently, the correspondence $F\mapsto \tilde{F}$ carries out an isomorphism between the categories $Shv(I, \Ab)$ and $\Ab^{I^{op}}$.
\hfill$\Box$

The sheaf generated by the presheaf $\mF$ takes the values 
$$
\hat\mF(U)=\invlim_{U^{op}} \{\mF(\Lambda_i)\}_{i\in U}.
$$
Lemma \ref{shpresh} allows one to identify sheaves over $I$ with diagrams of $I^{op}$.

\begin{proposition}\label{dlimsect}
$\invlim^n_{I^{op}} F \cong H^n(I, \widetilde{F})$.
\end{proposition}
{\sc Proof.}
The functor $F\mapsto \tilde{F}$ is exact. It takes the flasque resolution
$F^*$ of a diagram $F$ into the flasque resolution $\tilde{F}^*$ of a pencil 
$\tilde{F}$.
Since $\Gamma(\tilde{F})= \tilde{F}(I)= \invlim_{I^{op}}F$, the statement follows from the isomorphism of the complexes $\Gamma(\tilde{F^* })\cong
\invlim_{I^{op}}F^*$.
\hfill
$\Box$

\subsection{Examples of calculating {\v{C}ech} cohomology of posets}

Following \cite[III.3.8]{gro1957}, we define the {\v{C}ech} cohomology $\che^n(I,\mF)$ with coefficients in the presheaf $\mF: \tau_I^{op}\rightarrow \Ab$ as a colimit of $\dirlim H^n(\fU, \mF)$ over a directed family of classes of 
refinement open coverings.
Since there is the covering $\fU_I=\{\Lambda_i: i\in I\}$
refined all coverings, then 
$$ 
\che^n(I,\mF)=H^n(\fU_I,\mF), 
$$ 
and the study can be limited by considering only one covweing.
Consider the embedding of posets ${\Lambda^{op}}:I^{op}\xrightarrow{\subseteq} \widetilde{\fU}_I^{op}$.

Let us give an example showing that topos cohomology and {\v{C}ech} cohomology with coefficients in the sheaf for an Aleksandrov $T_0$-space may not be isomorphic.
\begin{example}\label{sphere2}
Consider the poset $S^{op}_2$:
$$
\xymatrix{
&
0 \ar[d] \ar[rrd] & & 1\ar[lld]\ar[d] \\
& 2\ar[d] \ar[rrd] && 3\ar[lld]\ar[d]\\
&4 && 5
}
$$
Let's take the sheaf $\mF(U)= \invlim_{U^{op}}F|_{U}$ on $S_2$ corresponding to the diagram $F= \Delta\ZZ$.
The geometric realization for $S_2$ is homeomorphic to the sphere $S^2$, whence $\invlim^2_{S^{op}_2}\Delta\ZZ\cong H^2(S^2, \ZZ)= \ZZ$.
By Theorem \ref{cheaslim} {\v{C}ech} cohomology theorem of the sheaf $\mF(U)= \invlim_{U^{op}} F|_U$ are equal to $\invlim^n_{\tilde\fU_{S_2}^{op}}\mF|_{\tilde\fU_{S_2}}$.
A poset whose elements are nonempty finite intersections is equal to $\tilde\fU_{S_2}= 
\{\Lambda_0, \Lambda_1, \Lambda_0\cap \Lambda_1,\Lambda_2, \Lambda_3, \Lambda_2\cap \Lambda_3, \Lambda_4, \Lambda_5\}$. It contains the subset $C_2= \{\Lambda_0, \Lambda_1, \Lambda_0\cap \Lambda_1\}$.
The full embedding of the category $C^{op}_2$ into ${\tilde\fU_{S_2}^{op}}$ is strictly coinitial, which implies $\check{H}^n(S_2, \mF)=0$, for all $n\geq 2$.
Consequently, the {\v{C}ech} cohomology of the space $S_2$ with coefficients in the
sheaf $\mF$ is not isomorphic to the topos cohomology of the sheaf $\mF$.
\end{example}

It has long been known \cite[III.3.8, Corollary 3]{gro1957} that for any sheaf
 $\mF$ on an arbitrary topological space $X$ there are isomorphisms
${\check H}^n(X, \mF) \to H^n(X,\mF)$ for $n=0$ and $n=1$. It is erroneous to say that this is true for presheaves.
The next two examples show that this is not the case when $\mF$ is a presheaf. Homomorphisms ${\check H}^n(X, \mF) \to H^n(X,\hat\mF)$ for $n=0$ and $n=1$ may not be isomorphisms. Here $\hat\mF$ denotes the sheaf generated by the presheaf $\mF$.
\begin{example}\label{contra01}
Let us calculate the {\v{C}ech} cohomology groups with coefficients in the presheaf over the partially ordered set $C_2$ from the example \ref{fourp}.
The set $\widetilde{\fU}_{C_2}^{op}$ will have the following Hasse diagram
$$
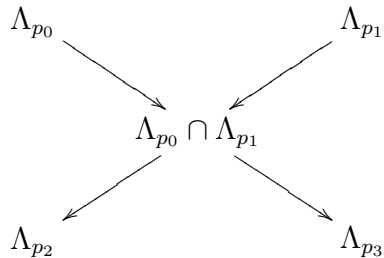
\xymatrix{
\Lambda_{p_0}\ar[rd] & & \Lambda_{p_1}\ar[ld]\\
& \Lambda_{p_0}\cap \Lambda_{p_1} \ar[ld] \ar[rd] \\
\Lambda_{p_2} & & \Lambda_{p_3}
}
$$
Consider, for example, a presheaf $\mF$ taking on $\widetilde{\fU}_{C_2}^{op}$ the values
$$
\xymatrix{
\ZZ\ar[rd]_{1_{\ZZ}} & & \ZZ\ar[ld]^{1_{\ZZ}}\\
& \ZZ \ar[ld] \ar[rd] \\
0 & & 0
}
$$
Since $\che^0(C_2,\mF)\cong \ZZ$ and $\invlim_{C_2}\mF\vert_{C_2}\cong \ZZ\oplus\ZZ$, then even the {\v{C}ech} $0$-cohomology groups of the presheaf 
$\mF$ are not isomorphic to the topos $0$-cohomology groups of the sheaf generated by the presheaf $\mF$.
In this case, the first cohomology groups are both equal to $0$.
Consider a presheaf $\mF$ taking the values
$$
\xymatrix{
\ZZ\ar[rd]_{1_{\ZZ}} & & \ZZ\ar[ld]^{1_{\ZZ}}\\
& \ZZ \ar[ld]_{1_{\ZZ}} \ar[rd]^{1_{\ZZ}} \\
\ZZ & & \ZZ
}
$$
Let $V\subset \widetilde{\fU}_{C_2}^{op}$ be a subset consisting of the three upper points.
It is strictly coinitial, hence $\che^1(C_2,\mF) \cong \invlim^1 \mF\vert_{V}= 0$.
The first group of topos cohomology according to Example \ref{fourp} will be isomorphic to the colimit of the diagram
$$
\xymatrix{
\ZZ \ar[d]_{1_{\ZZ}} \ar[rd]|-(.3){1_{\ZZ}} & \ZZ \ar[ld]|-(.3){1_{\ZZ}}\ar[d]^{1_{\ZZ}}\\
\ZZ & \ZZ
}
$$
Consequently, $\invlim^1 \mF\vert_{C_2} \cong \ZZ$, and the first {\v{C}ech} cohomology groups (with coefficients in the presheaf) and topos
cohomology are not isomorphic.
\end{example}

\subsection{Grothendieck and \v{C}ech cohomology
for Aleksandrov spaces}

Let $X\subseteq I$ be a subset of a poset.
Following \cite[IV, \S 11]{kur1968}, we denote
$$
X^{-}= \{i\in I: (\forall x\in X) i\leq x \},~
X^{+}= \{i\in I: (\forall x\in X) i\geq x \}.
$$
It is clear that $X^{-}= \bigcap\limits_{x\in X}\Lambda_x$.
This means that for each $W\in \widetilde{\fU}_I$ there is a finite set $S\subseteq I$ such that $W=S^-$.
If we denote $V_x= \{i\in I: i\geq x\}$, then we obtain a dual definition for $X^+= \bigcap\limits_{x\in X}V_x$.
Recall that $\hat\mF$ denotes the sheaf generated by the presheaf $\mF$.

\begin{theorem}\label{miso}
The following properties of an arbitrary poset $(I,\leq)$ are equivalent
\begin{itemize}
\item [1.] Canonical homomorphisms $\lambda_n: {\check H}^n(I, \mF)\to H^n(I, {\hat\mF})$ are isomorphisms for every presheaf $\mF$ Abelian groups over $I$ with the Aleksandrov topology and for each $n\geq 0$.
\item [2.]For any non-empty finite subset $X\subseteq I$ such that $X^-\not= \emptyset$, the integer homology groups $H_n(X^{-+})$ are isomorphic to the homology groups of the point for all $n\geq 0$.
\end{itemize}
\end{theorem}
{\sc Proof.}
From Theorem \ref{oberst} it follows that the following two conditions are equivalent:
\begin{itemize}
\item [1'.] Canonical homomorphisms $\invlim^n_{\tilde\fU^{op}}G \to \invlim^n_{I^{op}}G\Lambda^{op}$ associated with the functor $\Lambda: I\to \tilde\fU$, are isomorphisms for any diagram $G: \widetilde{\fU}^{op}\to \Ab$ and $n\geq 0$.
\item [2'.] For all $n\geq 0$ the integer homology groups $H_n(\Lambda/\Lambda_{x_1}\cap \cdots \cap \Lambda_{x_k}, \ZZ)$ are isomorphic to the point homology groups for any subset $\{x_1, \ldots, x_k\}\subseteq I$, $k\geq 1$, such that $\Lambda(x_1)\cap \cdots \cap \Lambda(x_k)\in \tilde\fU_I$.
\end{itemize}

 Let us prove the implication $1'\Rightarrow 1$.
To this end, let us substitute into condition $1'$ instead of the diagram $G: \widetilde{\fU}^{op}\to \Ab$ the restriction $\mF|_{\widetilde{\fU}}$ of an arbitrary presheaf $\mF$ over $(I, \tau_I)$.
If $1'$ is true, then the homomorphism
\begin{equation}\label{cond1}
\invlim^n_{\tilde\fU^{op}}\mF|_{\tilde\fU}\to \invlim^n_{I^{op}}
 \mF|_{\tilde\fU}\Lambda^{op}
\end{equation}
will be an isomorphism. The domain of this homomorphism will be isomorphic to $\che^n(I, \mF)$ by Theorem \ref{cheaslim}, and the codomain is isomorphic to $H^n(I, \hat\mF)$ by Proposition \ref{dlimsect}. So 1 is true.

Let condition 1 be satisfied.
Then we extend $G$ to the presheaf $\mF= Ran_{\tilde\fU^{op}\subseteq \tau^{op}_I}G$. For each $W\in \tilde\fU$ the comma category $W/\tilde\fU^{op}$ has an initial object, which implies that $\mF|_{\tilde\fU} =
(Ran_{\tilde\fU^{op}\subseteq \tau^{op}_I}G)|_{\tilde\fU}= G$.
  Since (\ref{cond1}) is an isomorphism, then 1' is true.

Let's prove $2\Leftrightarrow 2'$. For each finite set $X= \{x_1, \ldots, x_k\}$ the sets $X^-$ are equal to $\Lambda_{x_1}\cap \cdots \cap \Lambda_{x_k}$.
It follows that the sets ${X^-}^+$ are equal to the comma categories $\Lambda/\Lambda_{x_1}\cap \cdots \cap \Lambda_{x_k}$.
Condition 2' that the category is acyclic is equivalent to the fact that for each $X= \{x_1, \ldots, x_k\}$, such that the set $X^-=\Lambda(x_1)\cap \cdots \cap\Lambda(x_k )$ is not empty, the integer homology of the set $X^{-+}$ is isomorphic to the homology of the point.

We have proven $1\Leftrightarrow 1' $ and $2 \Leftrightarrow 2'$, from Theorem \ref{oberst} it follows $1'\Leftrightarrow 2'$. Therefore, $1\Leftrightarrow 2$.
\hfill
$\Box$
\\

According to \cite[IV, \S 11, (5)]{kur1968}, for each subset $X\subseteq I$ of a poset, the pair $\langle X^-, X^{-+}\rangle$ is a cut in $I$ and any cut can be represented in this form. This means that for finite posets $(I, \leq)$ the following statement is true.
\begin{corollary}\label{misof}
The following properties of a finite poset $I$ are equivalent
\begin{itemize}
\item Canonical homomorphisms $\lambda_n: {\check H}^n(I, \mF)\to H^n(I, {\hat\mF})$ are isomorphisms for every presheaf $\mF$ of Abelian groups over $ (I, \tau_I)$ and for each $n\geq 0$.
\item For any cut $\langle X, Y\rangle$ such that $X\not=\emptyset$, the integer homology groups $H_n(Y)$ are isomorphic to the homology groups of a point for all $n\geq 0$.
\end{itemize}
\end{corollary}

In particular, the condition of the theorem \ref{miso} is satisfied for poset directed from above, as well as for lower semilattices.
Let $I_0$ denote the partially ordered set obtained by adding the smallest element to $I$. If the smallest element exists in $I$, then we set $I_0=I$.

\begin{corollary}\label{semilat}
If each connected component $I$ is directed from above or $I_0$ is a lower semilattice, then the \v{C}ech cohomology of the topological space $(I,\tau_I)$ with coefficients in an arbitrary presheaf $\mF$ is isomorphic to the topos cohomology with coefficients in the sheaf generated by the presheaf $\mF$.
\end{corollary}
{\sc Proof.} In the first case, for every finite $X\subseteq I$ with non-empty $X^-$, the set $X^{-+}$ is non-empty and directed from above.
In the second case, for any $X=\{x_1, \cdots, x_n\}\subseteq I_0$ there is an infimum $x\in I_0$. True $\Lambda_{x_1}\cap \cdots \cap \Lambda_{x_n}=\Lambda_x$, and therefore $X^{-+}=\{x\}$.
If $x=0$, then $x\notin I$, whence $X^-=\emptyset$ in $I$.
This means that for every finite $X\subseteq I$ such that $X^-\not=\emptyset$, the second condition of Theorem \ref{miso} will be satisfied.
\hfill$\Box$

\subsection{Examples of comparison of \v{C}ech and topos cohomology of posets}

\begin{example}
Consider the partial ordered set $N$ defined using the Hasse diagram:
$$
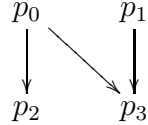
\xymatrix{
p_0 \ar[d] \ar[rd] & p_1\ar[d] \\
p_2 & p_3
}
$$
There is an isomorphism $\Lambda: N^{op}\stackrel{\cong}\rightarrow \widetilde{\fU}^{op}$.
Therefore, the \v{C}ech cohomology $\che^n(N,\mF)$ with coefficients in the arbitrary presheaf $\mF$ is isomorphic to the topos cohomology of the corresponding sheaf 
$\invlim^n \mF\vert_{N}$.
\end{example}

Let us give two examples of a finite poset for which the conditions of Theorem \ref{miso} are satisfied, but the conditions of the Corollary \ref{semilat} are not satisfied.

\begin{example}
Consider a poset $I$ defined using the Hasse diagram of the poset $I^{op}$
$$
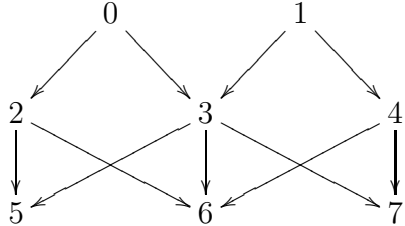
\xymatrix{
& 0\ar[ld] \ar[rd] & & 1\ar[ld]\ar[rd] \\
2\ar[d]\ar[rrd] & & 3\ar[lld]\ar[d]\ar[rrd] & & 4\ar[lld]\ar[d]\\
5 & & 6 & & 7
}
$$
We construct a covering $\tilde\fU_I$ from finite intersections of sets $\Lambda_i$ having incomparable indices.
It consists of the sets $\Lambda_i$, $0\leq i\leq 7$, as well as the intersections of pairs of the subset
$\Lambda_0 \cap \Lambda_1= \{3, 5, 6, 7\}$,
$\Lambda_0\cap \Lambda_4 = \Lambda_3\cap \Lambda_4 = \{6, 7\}$,
$\Lambda_1\cap \Lambda_2 = \Lambda_2\cap \Lambda_3 = \{5, 6\}$
and intersections of triples of subsets $\{2, 3, 4\}^-=\Lambda_2\cap \Lambda_3\cap \Lambda_4= \{6\}$ and $\{5, 6, 7\}^-= \Lambda_5 \cap\Lambda_6 \cap \Lambda_7 = \emptyset$.
Since $\{5,6\}^+= \Lambda/\{5,6\}= \{0, 1, 2, 3\}$ is a tree, it is connected and acyclic. Similar to $\{6,7\}^+$. In addition, $\{2,3,4\}^{-+}=\{3\}^+ = \{0, 1, 3\}$ is a tree.
From Theorem \ref{miso}, it follows that for each presheaf on the Aleksandrov space $(I, \tau_I)$ its \v{C}ech cohomology is isomorphic to the topos cohomology with the coefficients in the sheaf generated by this presheaf.
\end{example}

\begin{example} 
Consider a poset
$$
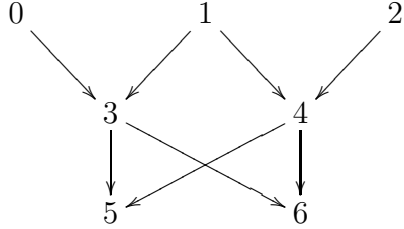
\xymatrix{
0 \ar[rd] && 1\ar[ld]\ar[rd] && 2\ar[ld]\\
& 3\ar[d] \ar[rrd] && 4\ar[lld]\ar[d]\\
&5 && 6
}
$$

Let us check the fulfillment of formula 2 from Theorem \ref{miso} for it.
For each $0\leq i\leq 6$, the poset $\{i\}^{-+}$ has the smallest element $i$, and hence it is acyclic.

Let's check for all pairs of incomparable elements $\{0,1\}^{-+}= \{0, 4\}^{-+}= \{0, 1, 3\}$, $\{0, 2 \}^{-+}= \{0, 1, 2, 3, 4\}$,
$\{1,2\}^{-+}= \{1, 2, 4\}$, $\{3,4\}^{-+}= \{0, 1, 2, 3, 4 \}$.
Let's check for a triple of elements $\{0,1, 2\}^{-+}= \{0, 1, 2, 3, 4\}$.
All tested upper sections are acyclic.
There is also a two-element subset such that $\{5, 6\}^-= \emptyset$.

For the constructed poset, \v{C}ech cohomology with coefficients in the presheaf and topos cohomology with coefficients in the sheaf generated by this presheaf are isomorphic.
\end{example}

\begin{example} 
Consider the cellular poset $P_X$ from \cite[Fig. 3]{eve2015}:
$$
\xymatrix{
& 0 \ar[ld]\ar[rd] \ar[rrrd] && 1\ar[llld] \ar[ld]\ar[rd] \\
2 \ar[d] \ar[rrd] && 3\ar[lld] \ar[rrd] && 4\ar[lld]\ar[d]\\
5\ar[rrd] && 6\ar[d] && 7\ar[lld]\\
&& 8
}
$$
Since $\{0,1\}^{-+}= \{0,1\}$, then the zero homology group $H_0(\{0,1\}^{-+})= \ZZ\oplus\ZZ$ is not isomorphic to $\ZZ$.
By Theorem \ref{miso} there exists $n\geq 0$ and a presheaf $\mF$ over the Aleksandrov space on $P_X$ for which the canonical homomorphism $\lambda_n: \che^n(P_{X}, \mF)\to H^n(P_X, {\hat\mF})$ will not be an isomorphism.
From \cite[Theorem 2]{eve2015} it follows that $H^n(P_X, {\hat{\mF}})\cong HC^n(P_X, {\hat{\mF}}\circ\Lambda^ {op})$,
where $HC^n(P, F)$ in \cite{eve2015} denotes the cellular cohomology of a partially ordered set $P$ with coefficients in the diagram $F: P^{op}\to \Ab$.
Therefore, the canonical mapping $\che^n(P_{X}, \mF)\to HC^n(P_X, {\hat\mF} \Lambda^{op})$ for these $n$ and $\mF$ is not an isomorphism.
\end{example}

\section{Conclusion}

The problem we posed in the introduction has been solved (Theorem \ref{miso}).

It remains unknown whether there is a criterion for the isomorphism of \v{C}ech cohomology for a sheaf on the $T_0$-Aleksandrov space and the topos cohomology of this sheaf, formulated in similar terms.
As Example \ref{sphere2} shows, in the general case the \v{C}ech cohomology and the topos cohomology of a sheaf may not be isomorphic.
The theorem \ref{miso} gives only sufficient conditions for isomorphism for sheaves.

In \cite{kis2022}, a theory of homology of co-sheaves on Aleksandrov spaces was constructed.
There is a problem to construct a theory of non-Abelian \v{C}ech homology for co-sheaves of groups on these spaces and to generalize for them sufficient conditions for the isomorphism of topos homology and \v{C}ech homology.
To solve it, you can suggest methods from the work \cite{X20233}.


\begin{thebibliography}{1}

\bibitem{bal2022}
V. Baladze, A. Beridze, L. Mdzinarishvili,
\textit{On axiomatic characterization of Alexander-Spanier normal homology theory of general topological spaces},
Topology Appl. 317, Article ID 108166, 25 p. (2022).


\bibitem{ben2020}
D. Bennequin, D. Peltre, G. Sergeant-Perthuis, J.P. Vigneaux,
\textit{Extra-fine sheaves and interaction decompositions},
New York, 2020. 47 p. (Preprint, Cornell Univ.) https://arXiv.org/abs/2009.12646

\bibitem{ber2020}
N. Berkouk, \textit{Persistence and Sheaves: from Theory to Applications}, Algebraic Topology [math.AT]. Institut Polytechnique de Paris, 2020.\\
https://theses.hal.science/tel-02970509v2/document


\bibitem{deh1962}
R. Deheuvels, \textit{Homologie des ensembles ordonn\'es
et des espaces topologiques}, Bulletin de la S. M. F., {\bf 90} (1962), 261--321.

\bibitem{eve2015}
B. Everitt, P. Turner, \textit{Cellular cohomology of posets with local coefficients},
Journal of Algebra, {\bf 439} (2015), 134--158.

\bibitem{gab1967}
P. Gabriel, M. Zisman,
\textit{Calculus of fractions and homotopy theory},
Berlin-Heidelberg-New York: Springer-Verlag, 1967.

\bibitem{god1958}
R. Godement,
\textit{Topologie alg\'ebrique et th\'eorie des faisceaux},
Herrman, Paris, 1958.

\bibitem{gro1957}
A. Grothendieck, \textit{Sur quelques points d'alg\`ebre homologique},
T\^ohoku Mathematical Journal, {\bf 9}: 2 (1957), 119-221.

\bibitem{X20232}
A. A. Husainov, \textit{Criteria for preserving the category cohomology
for the inverse image}, New York, 2023. 32 p. (Preprint, Cornell Univ.), \\
https://arXiv.org/abs/2306.12683

\bibitem{X20233}
A. A. Husainov, \textit{Homotopy cofinality for Non-Abelian homology of group diagrams}, New York, 2023. 25 p. (Preprint, Cornell Univ.), \\
https://arXiv.org/abs/2307.01726


\bibitem{jen1972}
C. U. Jensen,
     \textit{Les foncteurs d\'eriv\'es de
     $\invlim$ et leurs applications
     en th\'eorie des modules},
     Springer-Verlag, Berlin, 1972.

\bibitem{kas2018}
M. Kashivara, P. Schapira,
\textit{Persistent homology and microlocal sheaf theory.}
Journal of Applied and Computational Topology, {\bf 2} (2018), 83--113.

\bibitem{X1997}
A. A. Khusainov, \textit{Comparison of the dimensions of a small category},
Siberian Math. J. {\bf 38}:6 (1997), 1230--1240.

\bibitem{kis2022}
D. Kishimoto, Y. Yushima, \textit{Cellular cosheaf homology are cosheaf
homology}, New York, 2022. 11 p. (Preprint, Cornell Univ.);\\
https://arxiv.org/abs/2202.03659

\bibitem{kov1989}
V. A. Kovalevsky, \textit{Finite Topology as Applied to Image Analysis}, Computer vision, graphics, and image processing, {\bf 46} (1989),141--161.

\bibitem{kur1968}
K. Kuratowski, A. Mostowski,
\textit{Set theory},
PWN - Polish Scientific Publishers, Warszawa, 1968.

\bibitem{kuz1967}
V.I. Kuzminov,
\textit{Derived functors of the projective limit functor},
Sib. Mat. Zh., {\bf 8} (1967), 333--345. (Russian)

\bibitem{lau1965}
O.A. Laudal,  \textit{Sur les limites projectives et inductives},
Annales scientifiques de l'\'E.N.S. $\rm 3^e$ s\'erie, {\bf 82}: 2 (1965),
241--296.


\bibitem{mac2004}
S. MacLane,
\textit{Categories for the Working Mathematician},
Springer-Verlag, New York, 1998.

\bibitem{mir2011}
R. Miron, Gh. Piti\c{s}, I. Pop, \textit{On the Abstract \v{C}ech Cohomology},
Bul. Acad. \c{S}tiin\c{t}e Repub. Moldova. Matematica, {\bf 66}:2 (2011), 41--59.

\bibitem{obe1968}
U. Oberst,
        \textit{Homology of categories and exactness
        of direct limits},
        Math. Z., {\bf 107} (1968), 87--115.

\bibitem{san2020}
V.C. S\'anchez, C.M. P\'erez, F.S. de Salas, J.F.T. Sancho,
        \textit{Homology and cohomology of finite spaces},
 J. Pure Appl. Algebra, {\bf 224}:4 (2020), 106200.

\bibitem{sil2007}
V. de Silva, R. Ghrist, \textit{Coverage in sensor networks via persistent homology},
Algebraic and Geometric Topology, {\bf 7} (2007), 339-358.

\bibitem{tan2017}
K. Tanaka, \textit{Minimal networks for sensor counting problem using
discrete Euler calculus},
Japan J. Indust. Appl. Math., {\bf 34} (2017), 229-242.

\end{thebibliography}
\end{document}